\documentclass[a4paper,11pt]{amsart}

\usepackage{graphicx}
\usepackage{mathptmx}
\usepackage{amsmath}
\usepackage{amssymb}
\usepackage{enumitem}
\usepackage{xcolor}

\newmuskip\pFqmuskip

\newcommand*\pFq[6][8]{%
  \begingroup 
  \pFqmuskip=#1mu\relax
  \mathcode`=\string"8000
  \begingroup\lccode`\~=`\,
  \lowercase{\endgroup\let~}\pFqcomma
  F^{#2}_{#3}{\left(\genfrac..{0pt}{}{#4}{#5}\bigg|#6\right)}%
  \endgroup
}
\newcommand{\pFqcomma}{\mskip\pFqmuskip}

\newtheorem{theorem}{Theorem}[section]

\newtheorem{corollary}[theorem]{Corollary}
\newtheorem{proposition}[theorem]{Proposition}

\begin{document}

\title[Degenerate Eulerian polynomials and numbers]{Degenerate Eulerian polynomials and numbers}

\author{Dae San  Kim}
\address{Department of Mathematics, Sogang University, Seoul 121-742, Republic of Korea}
\email{dskim@sogang.ac.kr}
\author{Taekyun  Kim}
\address{Department of Mathematics, Kwangwoon University, Seoul 139-701, Republic of Korea}
\email{tkkim@kw.ac.kr}


\maketitle

\begin{abstract}
The aim of this paper is to study degenerate Eulerian polynomials and degenerate Eulerian numbers, respectively as degenerate versions of the Eulerian polynomials and the Eulerian numbers, and to derive some of their properties. Specifically, we derive an identity, recursive relations, generating function and degenerate version of Worpitzky's identity for the degenerate Eulerian polynomials and numbers. In addition, we obtain several results involving the degenerate Stirling numbers of the second kind and the degenerate Bernoulli numbers as well as the degenerate Eulerian numbers.
\end{abstract}

\section{Introduction}

The aim of this paper is to study the degenerate Eulerian polynomials $A_{n,\lambda}(x)$ (see \eqref{19}, \eqref{20}) and the degenerate Eulerian numbers $A_{n,\lambda}(n,k)$ (see \eqref{22}), respectively as a degenerate version of the Eulerian polynomials $A_{n}(x)$ (see \eqref{1}) and that of the Eulerian numbers $A_{n}(n,k)$ (see \eqref{3}). Among other things, we derive the properties in (a)-(e) below, all of which give alternative descriptions of the degenerate Eulerian polynomials or the degenerate Eulerian numbers. Moreover, we obtain several identities involving the degenerate Eulerian numbers, the degenerate Stirling numbers of the second kind and the degenerate Bernoulli numbers. \par
In more detail, the outline of this paper is as follows. In Section 1, we define the Eulerian polynomials $A_{n}(x)$ (see \eqref{1}) and the Eulerian numbers $A(n,k)$ (see \eqref{3}), which are the coefficients of $A_{n}(x)$, and recall an identity, recursive relations, generating function and Worpitzky's identity for the Eulerian polynomials and numbers, all giving alternative descriptions of those polynomials and numbers (see \eqref{5}-\eqref{9}). Then we remind the reader of the degenerate exponentials (see \eqref{11}) and the degenerate Bernoulli polynomials (see \eqref{12}). We recall the  degenerate Stirling numbers of the first kind (see \eqref{16}) and the second kind (see \eqref{18}).
Section 2 contains the main results of this paper. We define the degenerate Eulerian polynomials $A_{n,\lambda}(x)$ in \eqref{19} and the degenerate Eulerian numbers $A_{\lambda}(n,k)$ in \eqref{22}, respectively as a degenerate version of the Eulerian polynomials $A_{n}(x)$ and the Eulerian numbers $A(n,k)$. We derive (a)-(e) below, respectively in Proposition 2.1, Theorem 2.2, Theorem 2.3, theorem 2.6 and Theorem 2.7. They are respectively degenerate versions of the well-known properties on the Eulerian polynomials and numbers in \eqref{8}, \eqref{5}, \eqref{9}, \eqref{6} and \eqref{7}. In Theorem 2.4, we express $A_{n,\lambda}(-1)$ in terms of the degenerate Bernoulli numbers $\beta_{n,\lambda}$. In Theorems 2.8 and 2.9, we find expressions of the degenerate Stirling number of the second kind ${n \brace k}_{\lambda}$ and the sum $\sum_{k=1}^{m}(k)_{n,\lambda}$, both as finite sums involving $A_{-\lambda}(n,j)$. By combining $\sum_{k=1}^{m}(k)_{n,\lambda}=\frac{\beta_{n+1,\lambda}(m+1)-\beta_{n+1,\lambda}}{n+1}$ with the result in Theorem 2.9, we get the identity Theorem 2.10. Finally, we express $A_{\lambda}(n,k-1)$ as a finite sum involving ${n \brace j}_{\lambda}$.
\begin{align*}
&(a)\,\,\frac{x-1}{x-e_{-\lambda}((x-1)t)}=\sum_{n=0}^{\infty}A_{n,\lambda}(x)\frac{t^{n}}{n!}, \\
&(b)\,\,A_{\lambda}(n,k)=\sum_{i=0}^{k}\binom{n+1}{i}(-1)^{i}(k-i+1)_{n,\lambda},\,\,(0 \le k \le n),\\
&(c)\,\, A_{n,\lambda}(x)=\sum_{k=0}^{n-1}\binom{n}{k}A_{k,\lambda}(x)(1)_{n-k,-\lambda}(x-1)^{n-k-1}, \,\, (n \ge 1), \quad A_{0,\lambda}(x)=1,\\
&(d)\,\,A_{\lambda}(n,k)=\big((n-k)+(n-1)\lambda \big)A_{\lambda}(n-1,k-1)\\
&\qquad\qquad\qquad+\big(k+1-(n-1)\lambda \big))A_{\lambda}(n-1,k),\,\,(n \ge 1, \,0 \le k \le n),\\
&(e)\,\, (\textrm{Degenerate version of Worpitzky identity})\\
&\qquad \sum_{k=0}^{n}\binom{x+k}{n}A_{-\lambda}(n,k)=(x)_{n,\lambda},\,\,(n \ge 0).
\end{align*}

In combinatorics, the Eulerian number $A(n,k)$ is the number of permutations of the set $\{1,2,\dots,n\}$ in which exactly $k$ elements are greater than the previous element.
For example, there are 4 permutations of $\{1,2,3\}$ in which exactly 1 element is greater than previous element.\par
The Eulerian polynomials $A_{n}(x)$ were introduced in 1749 by Leonhard Euler in connection with computing values of the Riemann zeta function at negative integers (see [8]). The Eulerian polynomial $A_{n}(x)$ is the numerator of the rational function given by (see\ [7,11])
\begin{equation}
\sum_{j=0}^{\infty}(j+1)^{n}x^{j}=\frac{A_{n}(x)}{(1-x)^{n+1}},\quad (n \ge 0),\quad\mathrm{where}\ |x|<1.\label{1}
\end{equation}
We note that
\begin{equation}
\sum_{j=0}^{\infty}(j+1)^{n}x^{j}=\frac{1}{x}\bigg(x\frac{d}{dx}\bigg)^{n}\frac{1}{1-x},\quad (n \ge 1). \label{2}
\end{equation}
The Eulerian numbers $A(n,k)$ appear as the coefficients of the Eulerian polynomial $A_{n}(x)$. Indeed, we have
\begin{equation}
A_{n}(x)=\sum_{k=0}^{n}A(n,k)x^{k},\quad (\mathrm{see}\ [7,9,11]).\label{3}
\end{equation} \par
From \eqref{1}, we note that
\begin{align}
\sum_{k=0}^{n}A(n,k)x^{k}&=A_{n}(x)=\sum_{j=0}^{\infty}(j+1)^{n}x^{j}(1-x)^{n+1}\label{4}\\
&=\sum_{j=0}^{\infty}(j+1)^{n}x^{j}\sum_{l=0}^{\infty}\binom{n+1}{l}(-1)^{l}x^{l}	\nonumber\\
&=\sum_{k=0}^{\infty}\bigg(\sum_{i=0}^{k}\binom{n+1}{i}(-1)^{i}(k-i+1)^{n}\bigg)x^{k}. \nonumber
\end{align}
Thus, by comparing the coefficients on both sides of \eqref{4}, we get
\begin{equation}
\sum_{i=0}^{k}\binom{n+1}{i}(-1)^{i}(k-i+1)^{n}=\left\{\begin{array}{ccc}
A(n,k), & \textrm{if $0 \le k\le n$,}\\
0, & \textrm{if $k>n$,}
\end{array}\right.\quad (\mathrm{see}\ [7,11]). \label{5}
\end{equation} \par
The Eulerian numbers satisfy the recursive relation
\begin{equation}
A(n,k)=(n-k)A(n-1,k-1)+(k+1)A(n-1,k),\quad(n \ge 1, \,0 \le k \le n), \label{6}
\end{equation}
and \\
(Worpitzky's identity)
\begin{equation}
x^{n}=\sum_{k=0}^{n}A(n,k)\binom{x+k}{n},\quad (n \ge 0),\quad (\mathrm{see}\ [7,9,11]). \label{7}
\end{equation} \par
From \eqref{1}, we can derive the generating function of Eulerian polynomials given by
\begin{equation}
\frac{x-1}{x-e^{(x-1)t}}=\sum_{n=0}^{\infty}A_{n}(x)\frac{t^{n}}{n!},\quad (\mathrm{see}\ [7,11]).\label{8}
\end{equation}
By \eqref{8}, we easily get
\begin{equation}
A_{n}(x)=\sum_{k=0}^{n-1}\binom{n}{k}A_{k}(x)(x-1)^{n-1-k},\quad (n\ge 1),\quad A_{0}(x)=1,\label{9}
\end{equation}
and
\begin{equation}
\sum_{j=1}^{\infty}j^{n}x^{j}=\frac{xA_{n}(x)}{(1-x)^{n+1}}=\frac{1}{(1-x)^{n+1}}\sum_{k=0}^{n}A(n,k)x^{k+1},\quad (\mathrm{see}\ [5,7,11]). \label{10}	
\end{equation} \par
For any nonzero $\lambda\in\mathbb{R}$, the degenerate exponentials are defined by
\begin{equation}
e_{\lambda}^{x}(t)=\sum_{k=0}^{\infty}(x)_{k,\lambda}\frac{t^{k}}{k!},\quad e_{\lambda}^{1}(t)=e_{\lambda}(t),\quad (\mathrm{see}\ [11-20]), \label{11}
\end{equation}
where $(x)_{0,\lambda}=1,\ (x)_{k,\lambda}=x(x-\lambda)(x-2\lambda)\cdots\big(x-(k-1)\lambda\big),\ (k\ge 1)$. \\
Note that $\displaystyle \lim_{\lambda\rightarrow 0}e_{\lambda}^{x}(t)=e^{xt}\displaystyle$. \par
In [4], Carlitz considered the degenerate Bernoulli polynomials given by
\begin{equation}
\frac{t}{e_{\lambda}(t)-1}e_{\lambda}^{x}(t)=\sum_{n=0}^{\infty}\beta_{n,\lambda}(x)\frac{t^{n}}{n!}.\label{12}
\end{equation}
When $x=0$, $\beta_{n,\lambda}=\beta_{n,\lambda}(0),\ (n\ge 0),$ are called the degenerate Bernoulli numbers. The first few of the degenerate Bernoulli numbers are given by:
\begin{equation}
\beta_{0,\lambda}=1,\,\beta_{1,\lambda}=-\frac{1}{2}+\frac{1}{2}\lambda,\,  \beta_{2,\lambda}=\frac{1}{6}-\frac{1}{6}\lambda^{2},\, \beta_{3,\lambda}=-\frac{1}{4}\lambda+\frac{1}{4}\lambda^{3},\dots. \label{13}
\end{equation} \par
It is well known that the Stirling numbers of the first kind are defined by
\begin{equation}
(x)_{n}=\sum_{k=0}^{n}S_{1}(n,k)x^{k},\quad (n\ge 0),\quad (\mathrm{see}\ [1-7, 9-25]), \label{14}
\end{equation}
where $(x)_{0}=1,\ (x)_{n}=x(x-1)\cdots (x-n+1),\ (n\ge 1)$. \\
As the inversion formula of \eqref{14}, the Stirling numbers of the second kind are defined by
\begin{equation}
x^{n}=\sum_{k=0}^{n}{n \brace k}(x)_{k},\quad (n\ge  0),\quad (\mathrm{see}\ [1-7,9-25]).\label{15}
\end{equation} \par
Recently, the degenerate Stirling numbers of the first kind are defined by (see \eqref{14})
\begin{equation}
(x)_{n}=\sum_{k=0}^{n}S_{1,\lambda}(n,k)(x)_{k,\lambda},\quad (n\ge  0), (\mathrm{see}\ [13,14,15]). \label{16}
\end{equation}
As the inversion formula of \eqref{16}, the degenerate Stirling numbers of the second kind are defined by (see \eqref{15})
\begin{equation}
(x)_{n,\lambda}=\sum_{k=0}^{n}{n \brace k}_{\lambda}(x)_{k},\quad (n\ge  0),\quad (\mathrm{see}\ [13,14,15,18]).\label{17}
\end{equation}
An explicit expression of ${n \brace k}_{\lambda}$ is given by
\begin{equation}
{n \brace k}_{\lambda}=\frac{(-1)^k}{k!}\sum_{j=0}^{k}(-1)^{j}\binom{k}{j}(j)_{n,\lambda}. \label{18}
\end{equation}
Note that
\begin{displaymath}
\lim_{\lambda\rightarrow 0}S_{1,\lambda}(n,k)=S_{1}(n,k),\quad \lim_{\lambda\rightarrow 0}{n \brace k}_{\lambda}={n \brace k}.
\end{displaymath}

\section{Degenerate Eulerian polynomials and numbers}
In view of \eqref{1}, we define, as a degenerate version of $A_{n}(x)$, the {\it{degenerate Eulerian polynomials $A_{n,\lambda}(x)$}} by
\begin{equation}
\sum_{j=0}^{\infty}(j+1)_{n,\lambda}x^{j}=\frac{A_{n,\lambda}(x)}{(1-x)^{n+1}},\quad (n\ge 0),\quad\mathrm{where}\ |x|<1. \label{19}
\end{equation}
We observe that
\begin{equation}
\sum_{j=0}^{\infty}(j+1)_{n,\lambda}x^{j}=\frac{1}{x}\bigg(x\frac{d}{dx}\bigg)_{n,\lambda}\bigg(\frac{1}{1-x}\bigg),\quad (n \ge 1). \label{20}
\end{equation}
This is a degenerate version of \eqref{2}.
Note that $\displaystyle\lim_{\lambda\rightarrow 0}A_{n,\lambda}(x)=A_{n}(x),\ (n\ge 0)\displaystyle$. \par
Thus, by \eqref{19}, we get
\begin{align}
&\sum_{n=0}^{\infty}A_{n,\lambda}(x)\frac{t^{n}}{n!}=\sum_{n=0}^{\infty}\bigg(\sum_{j=0}^{\infty}(1+j)_{n,\lambda}(1-x)^{n+1}x^{j}\bigg)\frac{t^{n}}{n!}\label{21}\\
&=\sum_{j=0}^{\infty}x^{j}\sum_{n=0}^{\infty}(1+j)_{n,\lambda}(1-x)^{n+1}\frac{t^{n}}{n!}=(1-x)\sum_{j=0}^{\infty}x^{j}\sum_{n=0}^{\infty}(1+j)_{n,\lambda}(1-x)^{n}\frac{t^{n}}{n!}\nonumber \\
&=(1-x)\sum_{j=0}^{\infty}x^{j}e_{\lambda}^{j+1}\big((1-x)t\big)=(1-x)e_{\lambda}\big((1-x)t\big)\sum_{j=0}^{\infty}x^{j}e_{\lambda}^{j}\big((1-x)t\big)\nonumber\\
&=(1-x)e_{\lambda}\big((1-x)t\big)\frac{1}{1-xe_{\lambda}((1-x)t)}=\frac{1-x}{e_{\lambda}^{-1}((1-x)t)-x}=\frac{x-1}{x-e_{-\lambda}((x-1)t)}. \nonumber
\end{align}
Therefore, by \eqref{21}, we obtain the following theorem.
\begin{proposition}
\begin{displaymath}
\frac{x-1}{x-e_{-\lambda}((x-1)t)}=\sum_{n=0}^{\infty}A_{n,\lambda}(x)\frac{t^{n}}{n!}.
\end{displaymath}
\end{proposition}
\noindent This is a degenerate version of \eqref{8}. \par
In view of \eqref{3}, we define the {\it{degenerate Eulerian numbers}} $A_{n,\lambda}(n,k)$ by
\begin{equation}
A_{n,\lambda}(x)=\sum_{k=0}^{n}A_{\lambda}(n,k)x^{k},\quad (n\ge 0). \label{22}	
\end{equation}
Note that
\begin{displaymath}
\lim_{\lambda\rightarrow 0}A_{\lambda}(n,k)=A(n,k),\quad (n\ge k\ge 0).
\end{displaymath}
From \eqref{19} and \eqref{22}, we note that
\begin{align}
\sum_{k=0}^{n}A_{\lambda}(n,k)x^{k}=A_{n,\lambda}(x)&=\sum_{j=0}^{\infty}(j+1)_{n,\lambda}x^{j}(1-x)^{n+1}\label{23} \\
&=\sum_{j=0}^{\infty}(1+j)_{n,\lambda}x^{j}\sum_{i=0}^{\infty}\binom{n+1}{i}(-1)^{i}x^{i}\nonumber\\
&=\sum_{k=0}^{\infty}\bigg(\sum_{i=0}^{k}\binom{n+1}{i}(-1)^{i}(k-i+1)_{n,\lambda}\bigg)x^{k}.\nonumber
\end{align}
By comparing the coefficients on both sides of \eqref{23}, we get the following theorem.
\begin{theorem}
For $n,k\ge 0$, we have
\begin{displaymath}
\sum_{i=0}^{k}\binom{n+1}{i}(-1)^{i}(k-i+1)_{n,\lambda}=\left\{\begin{array}{ccc}
A_{\lambda}(n,k), & \textrm{if $0\le k\le n$,} \\
0, & \textrm{if $k>n$.}
\end{array}\right.
\end{displaymath}
\end{theorem}
\noindent This is a degenerate version of the identity in \eqref{5}. \par
From \eqref{21}, we note that
\begin{align}
x-1&=\sum_{i=0}^{\infty}A_{i,-\lambda}(x)\frac{t^{i}}{i!}\Big(x-e_{\lambda}\big((x-1)t\big)\Big)\label{24}\\
&=x\sum_{n=0}^{\infty}A_{n,-\lambda}(x)\frac{t^{n}}{n!}-\sum_{i=0}^{\infty}A_{i,-\lambda}(x)\frac{t^{i}}{i!}\sum_{l=0}^{\infty}(1)_{l,\lambda}(x-1)^{l}\frac{t^{l}}{l!}\nonumber\\
&=\sum_{n=0}^{\infty}\bigg(xA_{n,-\lambda}(x)-\sum_{i=0}^{n}\binom{n}{i}A_{i,-\lambda}(x)(1)_{n-i,\lambda}(x-1)^{n-i}\bigg)\frac{t^{n}}{n!}.\nonumber	
\end{align}
Thus, by \eqref{24}, we get
\begin{equation}
1=\sum_{n=0}^{\infty}\bigg(\frac{x}{x-1}A_{n,-\lambda}(x)-\sum_{i=0}^{n}\binom{n}{i}A_{i,-\lambda}(x)(1)_{n-i,\lambda}(x-1)^{n-i-1}\bigg)\frac{t^{n}}{n!}. \label{25}
\end{equation}
By comparing the coefficients on both sides of \eqref{25}, we have
\begin{equation}
\frac{x}{x-1}A_{n,-\lambda}(x)-\sum_{i=0}^{n}\binom{n}{i}A_{i,-\lambda}(x)(1)_{n-i,\lambda}(x-1)^{n-i-1}=\left\{\begin{array}{ccc}
1, & \textrm{if $n=0$,} \\
0, & \textrm{if $n>0$}.
\end{array}\right.	\label{26}
\end{equation}
Therefore, by \eqref{26}, we obtain the following theorem.
\begin{theorem}
For $n\ge 0$, we have
\begin{equation}
A_{0,\lambda}(x)=1,\quad A_{n,\lambda}(x)=\sum_{i=0}^{n-1}\binom{n}{i}A_{i,\lambda}(x)(1)_{n-i,-\lambda}(x-1)^{n-i-1},\quad (n\ge 1). \label{27}	
\end{equation}
\end{theorem}
\noindent This is a degenerate version of the recursive relation in \eqref{9}. \par
Letting $x=-1$ in \eqref{21} and recalling \eqref{13}, we have
\begin{align}
&\sum_{n=0}^{\infty}A_{n,\lambda}(-1)\frac{t^{n}}{n!}=\frac{-2}{-e_{-\lambda}(-2t)-1}=\frac{2}{1+e_{-\lambda}(-2t)}=\frac{2}{1+e_{\lambda}^{-1}(2t)} \label{28}\\
&=\frac{2}{e_{\lambda}(2t)+1}\big(e_{\lambda}(2t)+1-1\big)=2-\frac{2}{e_{\lambda}(2t)+1} \nonumber\\
&=2+\frac{4}{e_{\lambda}^{2}(2t)-1}-\frac{2}{e_{\lambda}(2t)-1}=2+\frac{4}{e_{\lambda/2}(4t)-1}-\frac{2}{e_{\lambda}(2t)-1}\nonumber \\
&=2+\frac{1}{t}\bigg(\frac{4t}{e_{\lambda/2}(4t)-1}-\frac{2t}{e_{\lambda}(2t)-1}\bigg)=2+\frac{1}{t}\bigg(\sum_{n=0}^{\infty}\Big(4^{n}\beta_{n,\lambda/2}-2^{n}\beta_{n,\lambda}\Big)\frac{t^{n}}{n!}\bigg)\nonumber\\
&=2+\frac{1}{t}\sum_{n=1}^{\infty}\Big(4^{n}\beta_{n,\lambda/2}-2^{n}\beta_{n,\lambda}\Big)\frac{t^{n}}{n!}=2+\sum_{n=0}^{\infty}2^{n+1}\bigg(\frac{2^{n+1}\beta_{n+1,\lambda/2}-\beta_{n+1,\lambda}}{n+1}\bigg)\frac{t^{n}}{n!}\nonumber \\ &=2+4\beta_{1,\lambda/2}-2\beta_{1,\lambda}+\sum_{n=1}^{\infty}2^{n+1}\bigg(\frac{2^{n+1}\beta_{n+1,\lambda/2}-\beta_{n+1,\lambda}}{n+1}\bigg)\frac{t^{n}}{n!}\nonumber\\
&=1+\sum_{n=1}^{\infty}2^{n+1}\bigg(\frac{2^{n+1}\beta_{n+1,\lambda/2}-\beta_{n+1,\lambda}}{n+1}\bigg)\frac{t^{n}}{n!}.\nonumber
\end{align}
Therefore, by comparing the coefficients on both sides of \eqref{28}, we obtain the following theorem.
\begin{theorem}
For $n\in\mathbb{N}$, we have
\begin{displaymath}
A_{0,\lambda}(-1)=1,\quad A_{n,\lambda}(-1)=2^{n+1}\bigg(\frac{2^{n+1}\beta_{n+1,\lambda/2}-\beta_{n+1,\lambda}}{n+1}\bigg).
\end{displaymath}
\end{theorem}
\begin{corollary}
For $n\in\mathbb{N}$, we have
\begin{displaymath}
\sum_{k=0}^{n}A_{\lambda}(n,k)(-1)^{k}=2^{n+1}\bigg(\frac{2^{n+1}\beta_{n+1,\lambda/2}-\beta_{n+1,\lambda}}{n+1}\bigg).
\end{displaymath}
\end{corollary}
From Theorem 2.2, we note that
\begin{align}
&A_{\lambda}(n,k)=\sum_{i=0}^{k}\binom{n+1}{i}(-1)^{i}(k-i+1)_{n,\lambda}\label{29}\\
&=\sum_{i=0}^{k}\binom{n+1}{i}(-1)^{i}(k-i+1)_{n-1,\lambda}\big(k-i+1-(n-1)\lambda\big) \nonumber\\
&=\big(k+1-(n-1)\lambda\big)\sum_{i=0}^{k}\binom{n+1}{i}(-1)^{i} (k-i+1)_{n-1,\lambda} \nonumber\\
&\qquad -\sum_{i=0}^{k}\binom{n+1}{i}(-1)^{i}i(k-i+1)_{n-1,\lambda}\nonumber
\end{align}
\begin{align*}
&=\big(k+1-(n-1)\lambda\big)\sum_{i=0}^{k}\bigg(\binom{n}{i-1}+\binom{n}{i}\bigg)(-1)^{i} (k-i+1)_{n-1,\lambda}\nonumber\\
&\qquad -(n+1)\sum_{i=1}^{k}\binom{n}{i-1}(-1)^{i} (k-i+1)_{n-1,\lambda}\nonumber\\
&=\big(k+1-(n-1)\lambda\big)\sum_{i=1}^{k}\binom{n}{i-1}(-1)^{i} (k-i+1)_{n-1,\lambda}\nonumber\\
&\qquad +\big(k+1-(n-1)\lambda\big)\sum_{i=0}^{k}\binom{n}{i}(-1)^{i}\big(k-i+1\big)_{n-1,\lambda} \nonumber\\
&\qquad +(n+1)\sum_{i=0}^{k-1}\binom{n}{i}(-1)^{i}\big(k-i\big)_{n-1,\lambda}\nonumber\\
&=-\big(k+1-(n-1)\lambda\big)A_{\lambda}(n-1,k-1)+\big(k+1-(n-1)\lambda\big)A_{\lambda}(n-1,k)\nonumber\\
&\qquad +(n+1)A_{k}(n-1,k-1). \nonumber
\end{align*}
Therefore, by \eqref{29}, we obtain the following theorem.
\begin{theorem}
For $n \ge 1$ and $0 \le k \le n$, we have
\begin{align}
&A_{\lambda}(n,k)=\big((n-k)+(n-1)\lambda\big)A_{\lambda}(n-1,k-1)\label{30}\\
&\qquad \qquad \quad+\big(k+1-(n-1)\lambda \big)A_{\lambda}(n-1,k).\nonumber
\end{align}
\end{theorem}
\noindent By taking $\lambda \rightarrow 0$, we get the recursive relation in \eqref{6}. \par
From \eqref{30}, we note that
\begin{align}
&\sum_{k=0}^{n+1}\binom{x+k}{n+1}A_{\lambda}(n+1,k)\label{31}\\
&=\sum_{k=0}^{n+1}\Big((n-k+1)A_{\lambda}(n,k-1)+(k+1)A_{\lambda}(n,k)\Big)\binom{x+k}{n+1}\nonumber\\
&\quad +n\lambda	\sum_{k=0}^{n+1}\Big(A_{\lambda}(n,k-1)-A_{\lambda}(n,k)\Big)\binom{x+k}{n+1}.\nonumber
\end{align}
Now, we observe that
\begin{align}
&\sum_{k=0}^{n+1}A_{\lambda}(n,k-1)\binom{x+k}{n+1}-\sum_{k=0}^{n+1}A_{\lambda}(n,k)\binom{x+k}{n+1}\label{32}\\
&=\sum_{k=1}^{n+1}A_{\lambda}(n,k-1)\binom{x+k}{n+1}-\sum_{k=0}^{n}A_{\lambda}(n,k)\binom{x+k}{n+1}\nonumber\\
&=\sum_{k=0}^{n}A_{\lambda}(n,k)\bigg(\binom{x+k+1}{n+1}-\binom{x+k}{n+1}\bigg)=\sum_{k=0}^{n}A_{\lambda}(n,k)\binom{x+k}{n},\nonumber	
\end{align}
and
\begin{align}
&\sum_{k=0}^{n+1}\Big((n-k+1)A_{\lambda}(n,k-1)+(k+1)A_{\lambda}(n,k)\Big)\binom{x+k}{n+1}\label{33}\\
&=\sum_{k=1}^{n+1}\binom{x+k}{n+1}(n-k+1)A_{\lambda}(n,k-1)+\sum_{k=0}^{n}A_{\lambda}(n,k)(k+1)\binom{x+k}{n+1}\nonumber\\
&=\sum_{k=0}^{n}\bigg((n-k)\binom{x+k+1}{n+1}A_{\lambda}(n,k)+\sum_{k=0}^{n}A_{\lambda}(n,k)(k+1)\binom{x+k}{n+1}\bigg)\nonumber\\
&=\sum_{k=0}^{n}A_{\lambda}(n,k)\bigg((n-k)\binom{x+k+1}{n+1}+(k+1)\binom{x+k}{n+1}\bigg). \nonumber
\end{align}
By simple calculation, we can show that
\begin{equation}
(n-k)\binom{x+k+1}{n+1}+(k+1)\binom{x+k}{n+1}
=x \binom{x+k}{n}.\label{34}
\end{equation}
Thus, by \eqref{33} and \eqref{34}, we get
\begin{equation}
\begin{aligned}
&\sum_{k=0}^{n+1}\Big((n-k+1)A_{\lambda}(n,k-1)+(k+1)A_{k}(n,k)\Big)\binom{x+k}{k}\\
&=x\sum_{k=0}^{n}A_{\lambda}(n,k)\binom{x+k}{n}.
\end{aligned}\label{35}
\end{equation}
From \eqref{31}, \eqref{32} and \eqref{35}, we have
\begin{equation}
\sum_{k=0}^{n+1}\binom{x+k}{n+1}A_{\lambda}(n+1,k)=(x+n\lambda)\sum_{k=0}^{n}\binom{x+k}{n}A_{\lambda}(n,k).\label{36}
\end{equation}
By applying \eqref{36} repeatedly, we get
\begin{align}
&\sum_{k=0}^{n+1}\binom{x+k}{n+1}A_{\lambda}(n+1,k) \label{37}\\
&=(x+n\lambda)\big(x+(n-1)\lambda\big)\cdots(x+\lambda)x\binom{x}{0}A_{\lambda}(0,0)=(x)_{n+1,-\lambda}.\nonumber
\end{align}
Therefore, by \eqref{37}, we obtain the following theorem.
\begin{theorem}[Degenerate version of Worpitzky's identity]
For $n\ge 0$, we have
\begin{equation*}
\sum_{k=0}^{n}\binom{x+k}{n}A_{-\lambda}(n,k)=(x)_{n,\lambda}.
\end{equation*}
\end{theorem}
\noindent By taking $\lambda \rightarrow 0$, we obtain the Worpitzky's identity in \eqref{7}. \par
From \eqref{17}, we note that
\begin{equation}
(x)_{n,\lambda}=\sum_{k=0}^{n}{n \brace k}_{\lambda}(x)_{k}=\sum_{k=0}^{n}k!{n \brace k}_{\lambda}\binom{x}{k},\quad (n\ge 0),\label{38}	
\end{equation}
and
\begin{equation}
\binom{x+y}{n}=\sum_{k=0}^{n}\binom{x}{k}\binom{y}{n-k},\quad (n\ge 0). \label{39}
\end{equation}
From Theorem 2.7, \eqref{38} and \eqref{39}, we have
\begin{align}
\sum_{k=0}^{n}k!{n \brace k}_{\lambda}\binom{x}{k}&=(x)_{n,\lambda}=\sum_{j=0}^{n}\binom{x+j}{n}A_{-\lambda}(n,j)\label{40}\\
&=\sum_{j=0}^{n}A_{-\lambda}(n,j)\sum_{k=0}^{n}\binom{x}{k}\binom{j}{n-k} \nonumber\\
&=\sum_{k=0}^{n}\binom{x}{k}\sum_{j=0}^{n}A_{-\lambda}(n,j)\binom{j}{n-k}. \nonumber	
\end{align}
Therefore, by comparing the coefficients on both sides of \eqref{40}, we obtain the following theorem.
\begin{theorem}
For $n\ge k\ge 0$, we have
\begin{displaymath}
{n \brace k}_{\lambda}=\frac{1}{k!}\sum_{j=0}^{n}A_{-\lambda}(n,j)\binom{j}{n-k}.
\end{displaymath}
\end{theorem}
For $m,n\in\mathbb{N}$ and Theorem 2.7, we have
\begin{align}
\sum_{k=1}^{m}(k)_{n,\lambda}&=\sum_{k=1}^{m}\sum_{j=0}^{n}\binom{k+j}{n}A_{-\lambda}(n,j)\label{41}\\
&=\sum_{j=0}^{n}	A_{-\lambda}(n,j)\sum_{k=1}^{m}\binom{k+j}{n}\nonumber\\
&=\sum_{j=0}^{n}A_{-\lambda}(n,j)\sum_{k=1}^{m}\bigg(\binom{k+j+1}{n+1}-\binom{k+j}{n+1}\bigg)\nonumber\\
&=\sum_{j=0}^{n}A_{-\lambda}(n,j)\binom{m+j+1}{n+1}. \nonumber
\end{align}
Therefore, by \eqref{41}, we obtain the following theorem.
\begin{theorem}
For $m,n\in\mathbb{N}$, we have
\begin{displaymath}
\sum_{k=1}^{m}(k)_{n,\lambda}=\sum_{j=0}^{n}A_{-\lambda}(n,j)\binom{m+j+1}{n+1}.
\end{displaymath}
\end{theorem}
From \eqref{12}, we note that
\begin{displaymath}
\beta_{n,\lambda}(x)=\sum_{k=0}^{n}\binom{n}{k}\beta_{k,\lambda}(x)_{n-k,\lambda},\quad (n\ge 0).
\end{displaymath}
Now, we observe that (see \eqref{12})
\begin{align}
\sum_{k=0}^{m}e_{\lambda}^{k}(t)&=\frac{e_{\lambda}^{m+1}(t)-1}{e_{\lambda}(t)-1}=\frac{1}{t}\bigg(\frac{t}{e_{\lambda}(t)-1}e_{\lambda}^{m+1}(t)-\frac{t}{e_{\lambda}(t)-1}\bigg)\label{42}\\
&=\frac{1}{t}\bigg(\sum_{n=0}^{\infty}\beta_{n,\lambda}(m+1)\frac{t^{n}}{n!}-\sum_{n=0}^{\infty}\beta_{n,\lambda}\frac{t^{n}}{n!}\bigg)\nonumber\\
&=\frac{1}{t}\sum_{n=1}^{\infty}\Big(\beta_{n,\lambda}(m+1)-\beta_{n,\lambda}\Big)\frac{t^{n}}{n!}\nonumber\\
&=\sum_{n=0}^{\infty}\bigg(\frac{\beta_{n+1,\lambda}(m+1)-\beta_{n+1,\lambda}}{n+1}\bigg)\frac{t^{n}}{n!}.\nonumber
\end{align}
From \eqref{42}, we get
\begin{equation}
\sum_{k=0}^{m}(k)_{n,\lambda}=\frac{\beta_{n+1,\lambda}(m+1)-\beta_{n+1,\lambda}}{n+1}, \label{43}
\end{equation}
where $m,n$ are non-negative integers.\\
Therefore, by \eqref{43} and Theorem 2.9, we obtain the following theorem.
\begin{theorem}
For $m,n\in\mathbb{N}$, we have
\begin{displaymath}
\sum_{j=0}^{n}A_{-\lambda}(n,j)\binom{m+j+1}{n+1}=\frac{1}{n+1}\Big(\beta_{n+1,\lambda}(m+1)-\beta_{n+1,\lambda}\Big).
\end{displaymath}
\end{theorem}
By using \eqref{18}, we observe that
\begin{align}
\bigg(x\frac{d}{dx}\bigg)_{n,\lambda}\bigg(\frac{1}{1-x}\bigg)&=\sum_{j=0}^{n}{n \brace j}_{\lambda}x^{j}\bigg(\frac{d}{dx}\bigg)^{j}\bigg(\frac{1}{1-x}\bigg)\label{44}\\
&=\sum_{j=0}^{n}\frac{(-1)^{j}}{j!}\sum_{k=0}^{j}(-1)^{k}\binom{j}{k}(k)_{n,\lambda}x^{j}\bigg(\frac{d}{dx}\bigg)^{j}\bigg(\frac{1}{1-x}\bigg)\nonumber\\
&=\sum_{j=0}^{n}(-1)^{j}\sum_{k=0}^{j}(-1)^{k}\binom{j}{k}(k)_{n,\lambda}\frac{x^{j}}{(1-x)^{j+1}},\nonumber
\end{align}
where $n$ is a positive integer. \\
From \eqref{44} and \eqref{20}, we have
\begin{equation}
\sum_{k=1}^{\infty}(k)_{n,\lambda}x^{k-1}=\frac{1}{x}\bigg(x\frac{d}{dx}\bigg)_{n,\lambda}\bigg(\frac{1}{1-x}\bigg)=\sum_{j=0}^{n}(-1)^{j}\sum_{k=0}^{j}(-1)^{k}\binom{j}{k}(k)_{n,\lambda}\frac{x^{j-1}}{(1-x)^{j+1}}.\label{45}
\end{equation}
By \eqref{19} and \eqref{45}, we get (see \eqref{10})
\begin{align}
\sum_{k=1}^{\infty}(k)_{n,\lambda}x^{k}&=\frac{xA_{n,\lambda}(x)}{(1-x)^{n+1}}=\frac{1}{(1-x)^{n+1}}\sum_{k=0}^{n-1}A_{\lambda}(n,k)x^{k+1}\label{46}\\
&=\frac{1}{(1-x)^{n+1}}\sum_{k=1}^{n}A_{\lambda}(n,k-1)x^{k}.\nonumber	
\end{align}
From \eqref{45} and \eqref{46}, we note that
\begin{align}
(1-x)\sum_{k=1}^{\infty}(k)_{n,\lambda}x^{k}&=\sum_{j=0}^{n}(-1)^{j}\sum_{k=0}^{j}(-1)^{k}\binom{j}{k}(k)_{n,\lambda}\frac{x^{j}}{(1-x)^{j}}\label{47} \\
&=\frac{1}{(1-x)^{n}}\sum_{k=1}^{n}A_{\lambda}(n,k-1)x^{k}.\nonumber	
\end{align}
Thus, by \eqref{47} and \eqref{18}, we get
\begin{align}
\sum_{k=0}^{n}A_{\lambda}(n,k-1)x^{k}&=\sum_{k=1}^{n}A_{\lambda}(n,k-1)x^{k}\label{48}\\
&=(1-x)^{n}\sum_{j=0}^{n}(-1)^{j}\sum_{k=0}^{j}(-1)^{k}\binom{j}{k}(k)_{n,\lambda}\frac{x^{j}}{(1-x)^{j}}\nonumber\\
&=(1-x)^{n}\sum_{k=0}^{n}k!{n \brace k}_{\lambda}\bigg(\frac{x}{1-x}\bigg)^{k}=\sum_{k=0}^{n}k!{n \brace k}_{\lambda}x^{k}(1-x)^{n-k}\nonumber\\
&=\sum_{k=0}^{n}k!{n \brace k}_{\lambda}x^{k}\sum_{j=0}^{n-k}\binom{n-k}{j}(-1)^{n-k-j}x^{n-k-j}\nonumber\\
&=\sum_{j=0}^{n}x^{n-j}\sum_{k=0}^{n-j}(-1)^{n-k-j}\binom{n-k}{j}k!{n \brace k}_{\lambda}\nonumber\\
&=\sum_{j=0}^{n}x^{j}(-1)^{j}\sum_{k=0}^{j}(-1)^{k}\binom{n-k}{n-j}k!{n \brace k}_{\lambda}\nonumber\\
&=\sum_{k=0}^{n}x^{k}(-1)^{k}\sum_{j=0}^{k}(-1)^{j}\binom{n-j}{n-k}j!{n \brace j}_{\lambda}. \nonumber
\end{align}
Therefore, by comparing the coefficients on both sides of \eqref{48}, we obtain the following theorem.
\begin{theorem}
For $n,k\in\mathbb{N}$, we have
\begin{displaymath}
A_{\lambda}(n,k-1)=(-1)^{k}\sum_{j=0}^{k}(-1)^{j}\binom{n-j}{n-k}j!{n \brace j}_{\lambda}.
\end{displaymath}
\end{theorem}
\section{Conclusion}
There are many ways of defining the Eulerian polynomials and numbers, including the ones in \eqref{1}, \eqref{3}, and \eqref{5}-\eqref{9}.
Another combinatorial interpretation of the Eulerian polynomial $A_{n}(x)$ is as follows. For each permutation $\sigma$ of the symmetric group $S_{n}$, the descent and the excedance of $\sigma$ are respectively given by
\begin{align*}
&d(\sigma)= |\{i \in [n-1] : \sigma(i) > \sigma(i + 1)\}|, \\
&e(\sigma)= |\{i \in [n-1] : \sigma(i) > i\}|,
\end{align*}
where $[n-1]=\{1,2,\dots,n-1\}$.
Then $A_{n}(x)$ is given by
\begin{align*}
A_{n}(x)=\sum_{\sigma \in S_{n}} x^{d(\sigma)}=\sum_{\sigma \in S_{n}} x^{e(\sigma)}.
\end{align*} \par
In this paper, after introducing the degenerate Eulerian polynomials and numbers we obtained degenerate version of the properties in (5)-(9), namely those in (a)-(e) in the Introduction, all giving alternative descriptions of those polynomials and numbers. Furthermore, we derived some results involving the degenerate Bernoulli numbers, the degenerate Stirling numbers of the second as well as the degenerate Eulerian numbers. \par
It is one of our future projects to continue to study degenerate versions, $\lambda$-analogues and probabilistic extensions of many special polynomials and numbers.

\end{document}